\documentclass[10pt, a4paper]{article}

\usepackage{amsmath,  amsfonts}

\usepackage[T1]{fontenc}

\newtheorem{Theorem}{\quad Theorem}[section]

\newtheorem{Corollary}[Theorem]{\quad Corollary}

\newtheorem{Example}[Theorem]{\quad Example}

\begin{document}

\title{ A  $ \sup + C \inf $ inequality on a domain of $ {\mathbb R}^2 $.}

 \author{Samy Skander Bahoura \footnote {e-mails: samybahoura@yahoo.fr, samybahoura@gmail.com} \\
 {\small Equipe d'Analyse Complexe et G\'eom\'etrie} \\  {\small Universit\'e Pierre et Marie Curie, 75005 Paris, France.}} 

\date{}

\maketitle

\begin{abstract}

Under some conditions, we give an inequality of type  $ \sup +  C\inf $ on open set of $ {\mathbb R}^2 $ for the prescribed scalar curvature equation.

\end{abstract}


\section{Introduction and Main Results} 

We set $ \Delta = \partial_{11} + \partial_{22} $  on open set $ \Omega $ of $ {\mathbb R}^2 $ with a smooth boundary.

Let's consider a sequence of functions solutions to:

$$ -\Delta u_i=V_i(x) e^{u_i}\,\,  \qquad (E) $$

with,

$$ 0 \leq V_i(x) \leq b,  \qquad (*)$$

$$ \int_{\Omega} e^{u_i} dx \leq C,  \qquad (**)$$

According to Brezis-Merle result, see [3], we have:

\smallskip

{\bf Theorem A}. {\it We have the following alternative:

1) $ (u_i)_i $ is uniformly locally bounded, 

\smallskip

or,

\smallskip

2) $ u_i  \to - \infty $ on each compact set of $ \Omega $, 

\smallskip

or,

\smallskip

3) There is a finite set (of blow-up points) $ S=\{ x_0, x_1, \ldots, x_m \} \subset \Omega $ and sequences of points $ (x_k^i)_i $ such that,

$$ x_k^i \to x_k, \,\,\, u_i(x_k^i) \to + \infty $$

and, in the sense of distributions we have:

$$ V_i e^{u_i} \to \sum_{k=0}^m \alpha_k \delta_{x_k},  \,\,\alpha_k \geq 4\pi, $$

and,

$$ u_i  \to - \infty , \,\, { \rm on \, each \, compact \, set \, of }\,\,  \Omega-S. $$

}

Here, we are interested by to know if,  we can have an inequality of type:

$$  \sup_{\Omega} u_i + C_1 \inf_K u_i \geq - C_2, $$

where $ C_1 $ and $ C_2 $ are two positive constant which depend only on $ b, C, K $ and $ \Omega $.

\smallskip
First, we give an example of  a corecif operator on a manifold with boundary for which the lower bound of the Green function on each compact set of the unit ball tends to 0.

\smallskip

\underbar { Nullity of the lower bound of the Green function for a manifold with boundary}

\smallskip

We have the following result which express that we can not apply the Nash-Moser iterate scheme to prove the existence of a lower bound for $ \sup + \inf  $ inequality. We need other arguments to obtain such estimate. Note that, in dimension $ n \geq  3 $, and, in the case of a compact riemannian manifold, we need the existence of a lower bound of the Green function to have a $ \sup \times \inf $ inequality. 

\smallskip

\begin{Example} Consider the sequence of blow-up functions :

$$ u_i (r) =\log \frac{8i^2}{(1+i^2r^2)^2} $$

which satisfy

$$ u_i(0) \to + \infty, \,\, { \rm and},\,\,  u_i(x) \to -\infty \,\, \forall \,\, x \not = 0, $$

and, for all  $ 0 <k \leq 1 $,

$$ u_i(0)+u_i(k) \to c_k \in (-\infty, +\infty), $$

Let's consider the Green function $ G_i $ of the following coercif operator in $ H_0^1(B_1(0)) $:

$$-\Delta +\epsilon_i, \,\, { \rm and},\,\,\epsilon_i=1+| \nabla u_i|^2,$$

with the following properties:

 $$  u_i'(r)=-\frac{4i^2r}{(1+i^2r^2)}, $$
 
 $$ || \nabla u_i||_{L^{\infty}(B_1(0)-B_r(0))} \leq C_r < +\infty, \,\, { \rm and},\,\,  | \nabla u_i(1/i)|=2i \to + \infty, $$

 then for all  $ 0 <k <1 $,

$$  \lim_i \left ( \inf_{x, y \in B_k(0)} G_i(x,y)\right ) = 0, $$

\end{Example}

We have the following theorem.
        
\begin{Theorem}  Let $ (u_i) $ and $ (V_i) $ two sequences of functions solutions to the previous problem, then

\bigskip

1) For all compact $ K  $ of $ \Omega $ we have:

$$ \sup_K u_i \to -\infty $$

or,

\smallskip

2) For all compact $ K \subset \Omega $, there are two positive constants $ C_1= C_1(dist(K, \partial \Omega), b, C) $ and $ C_2=C_2(K, \Omega, b, C) $ such that:

$$ \sup_{\Omega } u_i + C_1\inf_K u_i \geq -C_2 $$ 

\end{Theorem}

\section{Proof of the results:} 

\underbar {Proof of the Example 1: Nullity of the lower bound of the Green function.}

\bigskip

We want to prove by contradiction that: 

\bigskip

For all $ 1 > \epsilon >0 $ and $ \gamma >0 $ , there is a constant $ c= c(b, C, \epsilon, \gamma) $ such that:


$$ \gamma u_i(0) +  u_i(\epsilon)  \to + \infty, $$

For this, we consider the function:


 
 
 
$$   v_i=e^{u_i} $$ 
 
 $$   -\Delta  v_i + |\nabla u_i|^2 v_i=  v_i^2. $$

 Also, we can write:
 

$$   -\Delta v_i + (1+ |\nabla u_i|^2) v_i=  v_i^2+   v_i, $$

And,
 
  
  $$ \inf_{B_k(0)} v_i  \geq   \int_{B_{k+\epsilon}(0)} G_i(y_i,y)(v_i^2+  v_i ) dy $$

  Thus, by contradiction, we suppose that:  
   
   $$ \lim_i \left ( \inf_{x, y \in B_{k+\epsilon}(0)} G_i(x,y)\right )  > c >0, $$
   
   thus,   
   
  
  $$  \sup_{B_1(0)} e^{\gamma u_i}  \times v_i (k) \geq c  \int_{B_{k+\epsilon}(0)}(v_i^{2+\gamma}+ v_i^{1+\gamma} ) dy, $$

  If we argue by contradiction, and suppose that $ \sup + \inf  $ is finite:

  $$   \int_{B_{k+\epsilon}(0)}  v_i^{1+\gamma} dy   \leq C_1, \qquad (1) $$
  
  and,
  

$$ \int_{B_{k+\epsilon}(0)} v_i^{2+\gamma} dy \leq C_2, \qquad (2) $$      
  
  In particular, we have, using the Nash-Moser iterate scheme, the following inequality:

 $$ \int_{B_{k+\epsilon/2}(0)} |\nabla (\eta v_i^{\frac{1+\gamma}{2}})|^2 \leq C_3, \qquad (3) $$
     
   where $ \eta $ is a cutoff function.  With the Sobolev embeding, we obtain:
     
    $$ || v_i||_{L^q(B_{k+\epsilon/4})} \leq C_q \qquad (4) $$
 
 Also, we choose (because of Brezis-Merle result) $ r_0 $ such that, $ u_i \to  -\infty $ on $ \partial B_{r_0} $ and thus $ v_i \to 0 $ on $ \partial B_{r_0} $.
   
\bigskip 
   
 We use the Green representation formula to have:
  

    $$ v_i(0)=\int_{B_{r_0}(0)} G_0(0,y)(v_i^2-|\nabla  u_i|^2 v_i) dy + \int_{\partial B_{r_0}(0)} \partial_{\nu} G_0(0, s) v_i ds  $$

    Here, $ G_0 $ is the Green function of the Laplacian with Dirichlet condition on $  B_{r_0} $. We can write, $ G_0(x,y)=-\dfrac{1}{2\pi} \log |x-y|+H_0(x,y) $. For $ x \in  B_{1/2+\epsilon/6} $ and $ y \in  B_{r_0} $, $ G_0 \in L^q, \,\,\forall q \geq 1 $.  
   
    
    We use Holder inequality and $ (1) $, $ (2) $ or  $ 4 $ to obtain:    
       
    
    $$ v_i(0) \leq ||G_0(0,.)||_{L^q}||v_i||_{L^p}+ C_4||v_i||_{L^{\infty}(\partial B_{r_0})} \leq C_5 $$

 which contradict the fact that : 
 
 $$  v_i(0) =e^{u_i(0)} \to + \infty $$
 
 Thus, for all $ \gamma >0 $ :
 
    
    $$  \gamma u_i(0) + u_i(\epsilon) \to +\infty. $$    
    
 This means that, for our particular sequence:
 
 $$ u_i(0)+u_i(1) \to +\infty, $$  
 
 it is impossible. Thus:  
   
   $$ \lim_i \left ( \inf_{x, y \in B_{k+\epsilon}(0)} G_i(x,y)\right ) =0. $$
    
\underbar {Proof of the theorem 1}

\bigskip

Let's consider the sequence $ (u_i)_i $ on the ball of radius  2.

\bigskip

Now, suppose that, we are in the case 3) of the Brezis-Merle result. We assume that, $ \Omega= B_2(0) $, and, $ 0 \leq |x_0| \leq |x_1| \leq \ldots \leq |x_m| \leq 1/3 $  

We denote by $ G $ the Green function of the lapalcian on the unit ball.  We use the Green representation formula for $ u_i $, we obtain:

$$ u_i(y_i)=\max_{\Omega} u_i= \int_{B_1(0)} G(y_i,y) V_ie^{u_i(y)}dy +  \int_{\partial B_1(0)} \partial_{\nu} G(y_i, s) u_i ds  $$
 
 We can use the fact that (Brezis-Merle) $ u_i\to -\infty $ on $ \Omega-\cup_{k=1}^m B(x_k, r_k) $, to have:
 
 $$ \int_{B_1(0)} G(y_i,y) V_ie^{u_i(y)}dy=  \sum_{k=1}^m \int_{B(x_k, r_k)} G(y_i,y) V_ie^{u_i(y)}dy +\int_{\Omega-\cup_{k=1}^m B(x_k, r_k)} G(x_i,y) V_ie^{u_i(y)}dy = $$
 
 $$= \sum_{k=1}^m \int_{B(x_k, r_k)} G(y_i,y) V_ie^{u_i(y)}dy+ o(1), $$

Without loss of generality, we can assume that $ y_i \to x_0 $. Thus,  for $ k\not=0 $, $ G(y_i, y) \to \beta_k >0 $ for $ y \in B(x_k, r_k) $. We are concerned by the case $ y\in B(x_0, r_0) $. In this case, we can write:

$$   \int_{B(x_0, r_0)} G(y_i,y) V_ie^{u_i(y)}dy \leq \int_{B(y_i, r_0+ \epsilon)} G(y_i,y) V_ie^{u_i(y)}dy =\int_{B(y_i, r_0+ \epsilon)} \dfrac{1}{2\pi} \log \dfrac{|1-\bar y_iy|}{|y_i-y|} V_ie^{u_i(y)}dy $$
 
 But,
 
 $$ |1-\bar y_iy| \to \beta_0 >0, y \in B(x_0, r_0), $$
 
 Thus,
 
 $$   \int_{B(y_i, r_0+ \epsilon)} G(y_i,y) V_ie^{u_i(y)}dy \leq C(\beta, b, C)+ \int_{B(y_i, r_0+ \epsilon)} -\dfrac{1}{2\pi} \log |y_i-y| V_ie^{u_i(y)}dy $$ 
 
 We do a blow-up around $ y_i $, we write

  $$ y= y_i+xe^{-u_i(y_i)/2}, $$
  
  $$ \tilde u_i(x)= u_i( y_i+xe^{-u_i(y_i)/2})-u_i(y_i), $$ 
  
  and, 
  
  $$ \tilde V_i(x)= V_i( y_i+xe^{-u_i(y_i)/2}). $$
 
 We have,
 
 $$ \int_{B(y_i, r_0+ \epsilon)} -\dfrac{1}{2\pi} \log |y_i-y| V_ie^{u_i(y)}dy = \dfrac{u_i(y_i)}{4\pi} \int_{B(0, (r_0+ \epsilon) e^{u_i(y_i)/2})} \tilde V_ie^{\tilde u_i} dx + $$
 
 $$ + \int_{B(0, (r_0+ \epsilon) e^{u_i(y_i)/2})} -\dfrac{1}{2\pi} \log |x| \tilde V_i e^{\tilde u_i} dx , $$ 

The fact that $ e^{\tilde u_i} \leq 1 $ and $  -\log |x| \leq 0 $ for  $ |x| \geq 1 $, we can write:

  $$ \int_{B(0, (r_0+ \epsilon) e^{u_i(y_i)/2})} -\dfrac{1}{2\pi} \log |x| \tilde V_i e^{\tilde u_i} dx \leq C, $$ 
  
Also, we can see that, in the case $ V_i \to V  $ in $ C^0(\Omega) $,  by the result of YY.Li and I. Shafrir, we have:

$$ \int_{B(0, (r_0+ \epsilon) e^{u_i(y_i)/2})} \tilde V_ie^{\tilde u_i} dx =\int_{B(y_i, r_0+ \epsilon)}  V_ie^{u_i} dy \to 8\pi m \leq bC, \,\,m \in {\mathbb N}^* $$

Finaly, we can write:

$$ u_i(y_i) \leq (2m+\epsilon) u_i(y_i) + \sup_{\partial B_1(0)} u_i + C_1, $$

Thus, 

$$ C_2(bC) \sup_{\Omega} u_i +\sup_{\partial B_1(0)} u_i\geq -C_3, $$

But, 

$$ \sup_{\partial B_1(0)} u_i \to -\infty,  $$

Thus, by the classical Harnack inequality, we can have:

$$ \sup_{\partial B_1(0)} u_i \leq C_4 \inf_{\partial B_1(0)} u_i + C_5, $$

Finaly, for all compact set $ K $, there are two positive constants $ C_1'= C_1'(dist(K, \partial \Omega), b, C) $ and $ C_2'=C_2'(K, \Omega, b, C) $ such that:

$$ \sup_{\Omega} u_i + C_1'\inf_K u_i\geq -C_2', $$

In general, the condition $ V_i \to V  $ in $ C^0(\Omega) $ can be removed and replaced by the Brezis-Merle consequence of blow-up phenomenon.

$$ \int_{B(0, (r_0+ \epsilon) e^{u_i(y_i)/2})} \tilde V_ie^{\tilde u_i} dx =\int_{B(y_i, r_0+ \epsilon)}  V_ie^{u_i} dy \to 4\pi \tilde \alpha_0 \leq bC, \,\, \tilde \alpha_0 \geq 1. $$

\bigskip
            
    { \bf Remark}: {\it an example of the situation 1 of the theorem 2:}
     
     \bigskip
          
     Let's consider the following sequence of functions, on $ \Omega = B_1(0) $;
     
     $$ u_{\epsilon}(r)= \log \dfrac{1}{(\mu_{\epsilon}^2+ r^2)^2}  \,\,\, { \rm with} \,\,\mu_{\epsilon} \to +\infty $$  
          
     Then, for $ K = B_k(0) $ with $ 0 < k <1 $, we have:
     
     $$ -\Delta u_{\epsilon} = e^{u_{\epsilon}} \,\,\, { \rm and}, $$
     
     $$ \sup_{K} u_{\epsilon} \to - \infty, $$
     
     $$ \sup_{\Omega} u_{\epsilon} + \inf_K u_{\epsilon} \to - \infty  $$      
    
    { \bf Questions}:   1) a) Can we have the following sharp inequality:
    
    $$ \sup_{\Omega} u_i + \inf_K u_i \geq C_1=C_1(b, C, K, \Omega) \,\, ? $$    
   
 We look to the following example of I. Shafrir, see [7]:
 
   $$ u(r)= \begin{cases}

     2\log \left ( \dfrac{2\beta r^{\beta -1}}{1+r^{2\beta}}\right )\,\, \text{if } \,\, r >1 \,\, \\
     2\log \beta + 2\log \left ( \dfrac{2}{1+r^2}\right )\,\, \text{if } \,\, r \leq 1 .     
               
        \end{cases} $$

We take $ u_i(r)=u(ir)+2\log i $. Then, if we take $ \beta >1 $, we have:

$$ -\Delta u_i = \begin{cases}

     2 e^{u_i} \,\, \text{if }\,\, r >1/i \,\, \\
     \dfrac{2}{{\beta}^2}e^{u_i}\,\, \text{if } \,\, r \leq 1/i.     
          \end{cases}      $$
     
     $$  u_i(0) \to + \infty, $$
     
     $$ \forall \,\, 0 <k \leq 1,\,\, u_i(k) \to - \infty, $$
     
     $$ \forall \,\, 0 <k \leq 1,\,\, u_i(0) + u_i(k) \to - \infty. $$  
   
   b) Perhaps with more regularity on $ V_i $ ?
   
   With the previous example, we can use the example of Brezis-Li-Sahfrir, see [2],  to construct  two sequences of functions $ (u_i)_i $ and $ (V_i)_i $ ( $ 1< \beta=\beta_i \searrow 1 $),   such that:
   
   $$ -\Delta u_i = V_i e^{u_i} \,\, { \rm in} \,\, B_1(0), $$    
   
   $$ C^1(B_1(0)) \ni V_i \to V\equiv 2 \,\, { \rm in} \,\, C^0(B_1(0)), $$
   
   $$  u_i(0) \to + \infty, $$
     
     $$ \forall \,\, 0 <k \leq 1,\,\, u_i(k) \to - \infty, $$
    
    and we have:
     
     $$ \forall \,\, 0 <k \leq 1,\,\, u_i(0) + u_i(k) \to - \infty. $$  
      
   c) Perhaps, if we consider the solutions to the following equation:
   
   $$ -\Delta u_i = e^{u_i} \,\, { \rm in} \,\, \Omega $$  
     
    Consider the situation of X.Chen, see [5]. First we can have a blowing-up sequence $ \tilde u_i $, with two blow-up points, for example $ z_0=0 $ and $ z_1=1 $ and associeted radii $ r_i \to + \infty $, such that:

$$ e^{ \tilde u_i} \to \sum_{k=0}^1 8\pi \delta_{z_k}, , $$

and,

$$  \tilde u_i  \to - \infty , \,\, { \rm on \, each \, compact \, set \, of }\,\,  B_{r_i}-\{0, 1\}. $$

  $$  \int_{B_{r_i}} e^{ \tilde u_i} dx \to  16\pi. $$
  
  It is clear that we can choose the radius $ r_i $ such that :
  
  $$ r_i <<  \tilde u_i(0). $$
  
  Now, in the situation of the paper of X.Chen, we take the following sequence:
  
  $$ u_i(y)= \tilde u_i(r_iy)+ 2  \log r_i. $$
  
  In this case, we have one exterior blow-up point $ 0 $  and two interior blow-up points $ 0 $ and $ z_i=1/(r_i) $, in the unit ball, and, 
     
    $$ e^{u_i} \to 16\pi \delta_{0},  $$ 
     
     Let's consider the Green function of the unit ball, $ G $:
     
      $$ G(x, y)=\dfrac{1}{2 \pi} \log \dfrac{|1-\bar xy|}{|x-y|}, $$      
     
     We write:
     
     $$ u_i(0)=  \int_{B_1(0)} G(0,y) e^{u_i} dy +  \int_{\partial B_1(0)} \partial_{\nu} G(0, \sigma) u_i(\sigma) d\sigma, $$
     
     hence,
     
     $$ u_i(0)\geq  \int_{B_{l_i^0 r_i^0}(0)}  -\dfrac{1}{2 \pi} \log |y| e^{u_i} dy+ \int_{B_{l_i^1 r_i^1}(z_i)}  -\dfrac{1}{2 \pi} \log |y| e^{u_i} dy + \inf_{\partial B_1(0)} u_i, $$
     
         where,
         
         $$ r_i^0 = e^{-u_i(0)/2}, $$
         
         $$  r_i^1 = e^{-u_i(z_i)/2}, $$
         
         and, 
         
         $$ l_i^0 r_i^0 \leq  \dfrac{1}{2} |z_i|, $$
         
         $$ l_i^1 r_i^1  \leq   \dfrac{1}{2} |z_i|, $$
         
         we can choose $ l_i^0 $ as in the paper of CC.Chen and C.S Lin, see [4] (see also the formulation of Li-Shafrir, [6]), to obtain the following estimates:
         
         $$  \int_{B_{l_i^0 r_i^0}(0)}  -\dfrac{1}{2 \pi} \log |y| e^{u_i} dy = \dfrac{u_i(0)}{4 \pi} \int_{B_{l_i^0}}  e^{v_i} dt +  \int_{B_{l_i^0}}  -\dfrac{1}{2 \pi} \log |t| e^{v_i} dt,  $$         
         
         with,

            $$  \int_{B_{l_i^0}} e^{v_i} dt  \geq 8\pi-\dfrac{C}{L_i^s}-\dfrac{C}{u_i(0)}, $$
            
            and,
            
            $$ \int_{B_{l_i^0}}  \dfrac{1}{2 \pi} \log |t| e^{v_i} dt  \leq C,$$             
     
     where, $ v_i $ is the blow-up function around $ 0 $:
     
     $$ v_i(t)=u_i(e^{-u_i(0)/2} t)-u_i(0). $$
     
     With, 
     
     $$ (L_i)^s= ((1/2) |z_i| e^{u_i(0)/2})^s =((1/2 r_i) r_i e^{\tilde u_i(0)/2})^s=(1/2)^se^{s\tilde u_i(0)/2} >>  \tilde u_i(0)+ 2  \log r_i = u_i(0) . $$
     
     For the term:
     
     $$ \int_{B_{l_i^1 r_i^1}(z_i)}  -\dfrac{1}{2 \pi} \log |y| e^{u_i} dy,  $$    
     
     We take,
     
     $$ w_i(t)=u_i(e^{-u_i(z_i)/2} t+z_i)-u_i(z_i), $$   
      
      we have:
      
       $$ \int_{B_{l_i^1 r_i^1}(z_i)}  -\dfrac{1}{2 \pi} \log |y| e^{u_i} dy = \int_{B_{l_i^1}(0)} -\dfrac{1}{2 \pi} \log |t r_i^1+z_i|  e^{w_i} dt, $$
       
       which we can write as:
       
       $$ \int_{B_{l_i^1}(0)} -\dfrac{1}{2 \pi} \log |t r_i^1+z_i|  e^{w_i} dt =-\dfrac{1}{2 \pi} \log |z_i | \int_{B_{l_i^1}(0)} e^{w_i} dt+\int_{B_{l_i^1}(0)} -\dfrac{1}{2 \pi} \log |\dfrac{t r_i^1}{z_i}+1|  e^{w_i} dt, $$
       
       but, 
       
       $$ l_i^1 r_i^1  \leq   \dfrac{1}{2} |z_i|, $$
       
       Thus the term  $  \log |\dfrac{t r_i^1}{z_i}+1|  $ is bounded. And we have:

       $$ \int_{B_{l_i^1 r_i^1}(z_i)}  -\dfrac{1}{2 \pi} \log |y| e^{u_i} dy \geq  -\dfrac{1}{2 \pi} \log |z_i | \int_{B_{l_i^1}(0)} e^{w_i} dt- C, $$      
      
       because:
       
      $$  \int_{B_{l_i^1}(0)} e^{w_i} dt \to 8 \pi. $$
      
      Finally, we obtain the following inequality:
      
      $$ u_i(0) \geq 2 u_i(0)+  \inf_{\partial B_1(0)} u_i - ( 4-\epsilon) \log |z_i | - C, $$
      
      which we can write as:
      
      $$ u_i(0) +  \inf_{\partial B_1(0)} u_i \leq ( 4-\epsilon) \log |z_i | + C \to -\infty , $$

      Thus,
      
      $$ \sup_{\Omega } u_i + \inf_{ B_1(0)} u_i = u_i(0)+ \inf_{\partial B_1(0)} u_i \to -\infty. $$
       
       Here, we have choosen the Green function of the unit ball. The same argument may be applied to a ball of radius $ 1 \geq r > 0 $, to have:
       
     $$ \sup_{\Omega } u_i + \inf_{ B_r(0)} u_i = u_i(0)+ \inf_{\partial B_r(0)} u_i \to -\infty. $$         
           
   Thus, we have an example of blowing-up sequence with finite volume (conformal volume) and the $ \sup + \inf $ is not bounded below by a constant.      
   
 \bigskip
  
      2) What's happen if we remove  the condition $ (**) $ of the problem:
    
    $$   \int_{\Omega} e^{u_i} dx \leq C \,\,\, ?$$

Here, we can use the example:

$$ u(r)= \begin{cases}

     2\log \left ( \dfrac{2\beta r^{\beta -1}}{1+r^{2\beta}}\right )\,\, \text{if } \,\, r >1 \,\, \\
     2\log \beta + 2\log \left ( \dfrac{2}{1+r^2}\right )\,\, \text{if } \,\, r \leq 1 .     
               
        \end{cases} $$

We take $ u_i(r)=u(ir)+2\log i $. Then, if we replace $ \beta $ by $ i $, we have:

$$ -\Delta u_i = \begin{cases}

     2 e^{u_i} \,\, \text{if }\,\, r >1/i \,\, \\
     \dfrac{2}{i^2}e^{u_i}\,\, \text{if } \,\, r \leq 1/i.     
          \end{cases}      $$
     
     $$  u_i(0) \to + \infty, $$
     
     $$ \forall \,\, 0 <k \leq 1,\,\,  u_i(k) \to - \infty, $$
     
     $$ \forall \,\, C >0,  \forall \,\, 0 <k \leq 1,\,\, u_i(0) + C u_i(k) \to - \infty  $$ 

and,

 $$ \int_{B_1(0)} e^{u_i} dx \geq C' i \to + \infty, $$

3) Can we replace it, as in the paper of Brezis-Merle, by the following condition:

$$  0 < a \leq V_i(x) \leq b \,\,? \qquad (***) $$

Here, we can use the example:

$$ u_i(r)=2\log \left ( \dfrac{2i r^{i -1}}{1+r^{2i}}\right )\,\, \text{if } \,\, 1\leq r \leq 2. $$
  
We have:

$$ -\Delta u_i = 2 e^{u_i} \,\, \text{if }\,\, 1\leq r \leq 2. $$
     
     $$  u_i(1) \to + \infty, $$
     
     $$ \forall \,\, 1 <k \leq 2 \,\, u_i(k) \to - \infty, $$
     
     $$ \forall \,\, C >0,  \forall \,\, 1 <k \leq 2,\,\, u_i(1) + C u_i(k) \to - \infty  $$ 

and,

 $$ \int_{ \{ x, 1\leq |x| \leq 2\} } e^{u_i} dx \geq C' i \to + \infty, $$

Case when the $ \sup + \inf  $ inequality holds:
   
   \bigskip
   
  Consider on the unit ball, the radial solutions to the previous equation with the following condition:
 
   $$  u_i(0) \to + \infty, $$
   
   $$ \forall \,\, 0 <k \leq 1,\,\, u_i(k) \to - \infty, $$
   
   $$ \int_{B_1(0)} e^{u_i} dx \leq C. $$  
 
  Then we have:
  
  $$   \inf_{\partial B_1(0)} u_i =\sup_{\partial B_1(0)} u_i = u_i(1), $$
  
  This means:
  
  $$  \sup_{\partial B_1(0)} u_i-\inf_{\partial B_1(0)} u_i = 0, $$  
  
  All the conditions of the paper of YY.Li, see [7], (on riemannian surfaces) are satisfied and thus, we have:
  
  $$  \forall \,\, 0 <k \leq 1,\,\, |u_i(0)+ u_i(k)|  \leq c=c(C, k), $$
  
  For the general case when we assume that the prescribed scalar curvature $ V_i $ uniformly lipschitzian,  it is sufficient to prove the YY.Li condition:
  
  $$    \sup_{\partial B_1(0)} u_i-\inf_{\partial B_1(0)} u_i  \leq  B. $$
  
 Note that, in the case of compact riemannian surfaces without boundary,  the last condition hold and then, the $ \sup + \inf  $ inequality hold.  
 
 Now, assume we are on a compact riemannian surface without boundary. Assmue that the sequence blow-up and we have the following inequality:
 
 $$ \int_{M} V_i e^{u_i} dx \leq 8\pi. $$
 
 Then, we have a $ \sup + \inf  $ inequality. (For the 2-sphere, we have another proof of this fact).
 
 Now, assume that, on a compact surface without boundary,  we have:
 
  $$ \int_{M} V_i e^{u_i} dx \leq C, $$ 
  
  then,
  
  $$ C \sup_M u_i + \inf_M  u_i \geq -C_2=-C_2(b, C, M). $$
  
 If we use YY.Li condition, we can extend this result to the unit ball $ B_1(0) \subset {\mathbb R}^2 $. According to the proof of theorem 2, if we assume for a blowing-up sequence, that:
 
 $$  0 < a \leq V_i(x) \leq b, $$ 
 
 $$ \int_{B_1(0)} V_i e^{u_i} dx \leq 8\pi, $$ 
  
  and,
  
   $$    \sup_{\partial B_1(0)} u_i-\inf_{\partial B_1(0)} u_i  \leq  B. $$
  
  Then, the $ \sup + \inf  $ inequality hold.  We have exactly:
  
$$ \sup_{B_1(0)} u_i+\inf_{B_1(0)} u_i  \geq  -C=-C(a, b, B). $$

The last case correspond to the case of one blow-up point.
 
 \bigskip
 
4) Is it possible to have $ C_1 $ independant of the compact $ K $ ? This is the case if we suppose the YY.Li condition satisfied.

\bigskip
 
For example, on a  compact riemannian surface $ M $ without boundary,  if we consider the following equation:

$$  -\Delta u_i + R = V_ie^{u_i}, $$

with, $ V_i $ such that:

$$  0 < a \leq V_i(x) \leq b, \,\, { \rm on} \,\, M, $$

we have after integration of the equation that:

$$ \forall \,\, \Omega \subset M, \,\,  0 < a \int_{\Omega} e^{u_i} dx \leq a \int_{M} e^{u_i} dx \,\, \leq \int_{M} R(x) dx \leq |M| \sup_{M} R, $$

Also, we have:

$$ 0 < \int_{M} R(x) dx = \int_{M} V_i e^{u_i} dx \leq b e^{\sup_{M} u_i}, $$

thus,

$$ \sup_{M} u_i \geq \log \left ( \dfrac{\int_{M} R(x) dx}{b} \right ) . $$

Thus, we are in the case of theorem 2 and the point  1) of this theorem is not possible. We have the following:

\begin{Corollary}  Let $ (u_i) $ and $ (V_i) $ two sequences of functions solutions to the previous problem on a riemannian surface $ M $ , such that, we have only the following condition:

$$  0 < a \leq V_i(x) \leq b,  \qquad (***) $$

 then, there are two positive constants $ C_1=C_1(a, b, M) $ and $ C_2=C_2(a, b, M) $ such that:

$$ \sup_{M} u_i + C_1\inf_M u_i \geq -C_2, $$ 

\end{Corollary}

 { \bf Question}: Note that in the previous paper, by using the maximum principle,  we have proved that if we consider on compact riemannian surface without boundary the following equation:
 
 $$   -\Delta u_i + R = V_ie^{u_i}, $$

with, $ V_i $ such that:

$$   a=0 \leq V_i(x) \leq b, \,\, { \rm on} \,\, M, $$

Then we have an inequality of type:

$$ \sup_{M} u_i + C_1\inf_M u_i \geq -C_2, $$

with $ C_1=C_1(b, M) $ and $ C_2=C_2(b, M) $.

\bigskip

Can we try to find this result with the previous argument of the theorem 2 ?
 
 We have,
 
 $$  \int_{M} V_i e^{u_i} dx  = \int_{M} R(x) dx \leq |M| \sup_{M} R = C, $$
 
 We can say that, after using the proof of theorem 2 and the  YY.Li inequality around the maximum of $ u_i $:
 
  $$ \forall \,\Omega, \,\, \sup_{\partial \Omega} u_i-\inf_{\partial \Omega} u_i  \leq  C',  $$ 

that,

\begin{Corollary}  Let $ (u_i) $ and $ (V_i) $ two sequences of functions solutions to the previous problem on a riemannian surface $ M $ without boundary, such that, we have only the following condition:

$$  0  \leq V_i(x) \leq b,  \qquad (****) $$

 then, there is a positive constant $ C_1=C_1(b, M) $ such that:

$$ \left (\dfrac{\int_M R}{4\pi}-1 \right ) \sup_{M} u_i + \inf_M u_i \geq -C_1, $$ 

\end{Corollary}

This is the wellknown inequality in a previous paper. To see this, we have assumed in a previous paper that:

$$ R \equiv k \in {\mathbb R}_+^*, \,\, {\rm and } \,\, Volume (M)=1, $$

thus, the inequality is :

$$ \left (\dfrac{k-4\pi}{4\pi} \right ) \sup_{M} u_i + \inf_M u_i \geq -C_1. $$

















\bigskip

\begin{thebibliography}{99} 

\bibitem{1}{T. Aubin. Some Nonlinear Problems in Riemannian Geometry. Springer-Verlag 1998 }

\bibitem{2}{H. Brezis, YY. Li , I. Shafrir. A sup+inf inequality for some
nonlinear elliptic equations involving exponential
nonlinearities. J.Funct.Anal.115 (1993) 344-358.
}

\bibitem{3}{H.Brezis and F.Merle, Uniform estimates and blow-up bihavior for solutions of $ -\Delta u=Ve^u $ in two dimensions, Commun Partial Differential Equations 16 (1991), 1223-1253.
}
\bibitem{4}{C-C.Chen, C-S. Lin. A sharp sup+inf inequality for a nonlinear elliptic equation in ${\mathbb R}^2$.
Commun. Anal. Geom. 6, No.1, 1-19 (1998).}
\bibitem{5}{X.Chen, Remarks on the existence of branch bubbles on the blowup analysis of equation $-\Delta u=e^{2u} $ in dimension two. Comm. Anal. Geom. 7 (1999), no. 2, 295Ð302.      }

\bibitem{6}{YY. Li, I. Shafrir. Blow-up Analysis for Solutions of $
-\Delta u = V e^u $ in Dimension Two. Indiana. Math. J. Vol 3, no 4.
(1994). 1255-1270.}

\bibitem{7}{YY. Li. Harnack Type Inequality: the Method of Moving Planes. Commun. Math. Phys. 200,421-444 (1999).}

\bibitem{8}{I. Shafrir. A sup+inf inequality for the equation $ -\Delta u=Ve^u $. C. R. Acad.Sci. Paris S\'er. I Math. 315 (1992), no. 2, 159-164.}


\end{thebibliography}
\end{document}